\documentclass[11pt,reqno]{amsart}

\usepackage{graphicx}
\usepackage[top=20mm, bottom=20mm, left=30mm, right=30mm]{geometry}
\usepackage[colorlinks=true,citecolor=black,linkcolor=black,urlcolor=magenta]{hyperref}

\usepackage{fancyvrb}
\usepackage{ulem} 
\usepackage[T1]{fontenc} 
\usepackage{pbsi} 
\usepackage{url} 

\usepackage{local-math-commands}

\usepackage{diagbox}
\newcommand{\trianglenk}[2]{$\diagbox{#1}{#2}$}

\newcommand{\HZCf}[4]{\gkpSI{#1}{#2}_{\left(#3, #4\right)}}
\newcommand{\HZCfStar}[4]{\gkpSII{#1}{#2}_{\left(#3, #4\right)^\ast}}
\newcommand{\HZNumber}[4]{H_{#1}^{\left(#2\right)}\left(#3, #4\right)}
\newcommand{\HZFact}[3]{#1!_{\left(#2, #3\right)}}
\newcommand{\SIIStarf}[3]{\gkpSII{#1}{#2}_{#3^\ast}}

\DeclareMathOperator{\csch}{csch} 

\title[Generating Function Transformations]{
       Zeta Series Generating Function Transformations Related to 
       Generalized Stirling Numbers and 
       Partial Sums of the Hurwitz Zeta Function}  

\author{Maxie D. Schmidt \\ 
        \href{mailto:maxieds@gmail.com}{maxieds@gmail.com}}         
\date{2016.06.15-v1} 

\allowdisplaybreaks

\begin{document}

\maketitle

\begin{abstract} 
We define a generalized class of modified zeta series transformations 
generating the partial sums of the Hurwitz zeta function and 
series expansions of the Lerch transcendent function. 
The new transformation coefficients we define within the article 
satisfy expansions by generalized harmonic number sequences, or the 
partial sums of the Hurwitz zeta function, which are 
analogous to known properties 
for the Stirling numbers of the first kind and for the 
known transformation coefficients employed to enumerate variants of the 
polylogarithm function series. 
Applications of the new results we prove in the article include 
new series expansions of the Dirichlet beta function, the 
Legendre chi function, BBP--type series identities for special constants, 
alternating and exotic Euler sum variants, alternating zeta functions 
with powers of quadratic denominators, and particular series 
defining special cases of the Riemann zeta function constants at the 
positive integers $s \geq 3$. 
\end{abstract}

\section{Introduction}
\label{Section_Intro} 

\subsection{Definitions} 

The \emph{generalized $r$--order harmonic numbers}, 
$\HZNumber{n}{r}{\alpha}{\beta}$, are defined as the partial sums of the 
\emph{modified Hurwitz zeta function}, $\zeta(s, \alpha, \beta)$, 
defined by the series 
\begin{align} 
\label{eqn_GenModifiedHZetaFn_def} 
\zeta(s, \alpha, \beta) & = \sum_{k \geq 1} \frac{1}{(\alpha k + \beta)^s}. 
\end{align} 
That is, we define these generalized sequences as the sums 
\begin{align} 
\label{eqn_GenHZNums_def} 
\HZNumber{n}{r}{\alpha}{\beta} & = \sum_{1 \leq k \leq n} 
     \frac{1}{(\alpha k+\beta)^r}, 
\end{align} 
where the definition of the 
``ordinary'' \emph{$r$--order harmonic numbers}, 
$H_n^{(r)}$, is given by the special cases of \eqref{eqn_GenHZNums_def} 
where $H_n^{(r)} \equiv \HZNumber{n}{r}{1}{0}$ \citep[\S 6.3]{GKP}. 
Additionally, we define the analogous ``\emph{modified}'' 
\emph{Lerch transcendent function}, 
$\Phi(z, s, \alpha, \beta) \equiv \alpha^{-s} \cdot \Phi(z, s, \beta / \alpha)$, 
for $|z| < 1$ or when $z \equiv -1$ by the series 
\begin{align} 
\label{eqn_ModLTransFn_series_def} 
\Phi(z, s, \alpha, \beta) & = \sum_{n \geq 0} 
    \frac{z^n}{(\alpha n+\beta)^s}. 
\end{align} 
We notice the particular important interpretation that the 
Lerch transcendent function 
acts as an ordinary generating function that enumerates the 
generalized harmonic numbers in \eqref{eqn_GenHZNums_def} 
according to the coefficient identity 
\begin{align*} 
\HZNumber{n}{r}{\alpha}{\beta} & = [z^n] 
     \frac{\Phi(z, s, \alpha, \beta) - \beta^{-s}}{1-z},\ 
     |z| < 1 \vee z = -1,\ n \geq 0. 
\end{align*} 

\subsection{Approach to Generating the Modified Zeta Function Series} 

The approach to enumerating the harmonic number sequences and 
series for special constants within this article begins in 
Section \ref{subSection_GenHZStirlingNumbers} 
with a brief overview of the properties of the harmonic number expansions in 
\eqref{eqn_GenHZNums_def} obtained through the definition of a 
generalized Stirling number triangle extending the results in 
\citep{MULTIFACTJIS}. 
We can also employ transformations of the generating functions of many 
sequences, though primarily of the geometric series, to enumerate and 
approximate the generalized harmonic number sequences in 
\eqref{eqn_GenHZNums_def} which form the partial sums of the 
modified Hurwitz zeta function in \eqref{eqn_GenModifiedHZetaFn_def}. 

For example, given the ordinary generating function, $A(t)$, of the 
sequence $\langle a_n \rangle_{n \geq 0}$, we can employ a known 
integral transformation involving $A(t)$ termwise to enumerate the following 
modified forms of $a_n$ when $r \geq 2$ is integer-valued 
\citep{EXPLICIT-EVAL-ESUMS}: 
\begin{align*}
\sum_{n \geq 0} \frac{a_n}{(n+1)^r} z^n & = 
     \frac{(-1)^{r-1}}{(r-1)!} \int_0^1 \log^{r-1}(t) A(tz) dt. 
\end{align*} 
For integers $\alpha \geq 2$ and $0 \leq \beta < \alpha$, 
we can similarly transform the 
ordinary generating function, $A(t)$, through the previous integral 
transformation and the \emph{$\alpha^{th}$ primitive root of unity}, 
$\omega_{\alpha} = \exp(2\pi\imath / \alpha)$, 
to reach an integral transformation for the 
modified Lerch transcendent function in \eqref{eqn_ModLTransFn_series_def} 
in the form of the next equation \citep{GFOLOGY}. 
\begin{align*} 
\sum_{n \geq 0} \frac{a_{\alpha n+\beta}}{(\alpha n + \beta +1)^r} 
     z^{\alpha n+\beta} & = 
     \frac{(-1)^{r-1}}{\alpha \cdot (r-1)!} \int_0^1 \log^{r-1}(t) \left( 
     \sum_{0 \leq m < \alpha} \omega_\alpha^{-m\beta} 
     A\left(\omega_{\alpha}^m tz\right) 
     \right) dt 
\end{align*} 
The main focus of this article is on the applications and expansions of 
generating function transformations introduced in 
Section \ref{subSection_Intro_GenGFSeriesEnums} 
that generalize the forms of the coefficients defined in \citep{GFTRANS2016}. 
In many respects, the generalizations we employ 
here to transform geometric--series--based generating functions into 
series in the form of \eqref{eqn_ModLTransFn_series_def} are 
corollaries to the results in the first article. 
The next examples illustrate the new series we are able to obtain using 
these generalized forms of the generating function transformations 
proved in the reference. 

\subsection{Examples of the New Results} 
\label{subSection_Intro_Examples} 

\subsubsection{BBP-Type Formulas and Identities} 

Many special constants such as those given as examples in the next 
sections satisfy series expansions given by \emph{BBP--type formulas} 
of the form \citep{BBP-FORMULAS} 
\begin{equation*} 
P(s, b, m, A) = \sum_{k \geq 0} b^{-k} \sum_{1 \leq j \leq m} 
     \frac{a_j}{(mk+j)^s} 
\end{equation*}
where $A = (a_1, a_2, \ldots, a_m)$ is a vector of constants and 
$m, b, s \in \mathbb{Z}^{+}$. 
A pair of particular examples of first-order BBP-type formulas 
which we list in this section to demonstrate the 
new forms of the generalized coefficients listed in 
Table \ref{table_TableOfGenCoeffs_ab31} and 
Table \ref{table_TableOfGenCoeffs_ab32} of the article below 
provide series representations for a real--valued multiple of $\pi$ and 
a special expansion of the natural logarithm function 
\cite[\S 11]{BBP-FORMULAS} \citep[\S 3]{DISC-BBPTYPE-FORMULAS}. 
\begin{align*}
\frac{4\sqrt{3} \pi}{9} & = 
     \sum_{k \geq 0} \left(-\frac{1}{8}\right)^k \left( 
     \frac{2}{(3k+1)} + \frac{1}{(3k+2)} \right) \\ 
     & = 
     \sum_{j \geq 0} \frac{8}{9^{j+1}} \left( 
     2 \binom{j+\frac{1}{3}}{\frac{1}{3}}^{-1} + \frac{1}{2} 
     \binom{j+\frac{2}{3}}{\frac{2}{3}}^{-1} \right) \\ 
\Log\left(\frac{n^2-n+1}{n^2}\right) & = 
     \sum_{k \geq 0} \left(-\frac{1}{n^3}\right)^{k+1} \left[ 
     \frac{n^2}{3k+1} - \frac{n}{3k+2} - \frac{2}{3k+3} \right] \\ 
     & = 
     - \sum_{j \geq 0} \frac{1}{(n^3+1)^{j+1}} \left( 
     n^2 \binom{j+\frac{1}{3}}{\frac{1}{3}}^{-1} - 
     \frac{n}{2} \binom{j+\frac{2}{3}}{\frac{2}{3}}^{-1} - 
     \frac{2}{3 (j+1)} \right) 
\end{align*} 

\subsubsection{New Series for Special Zeta Functions} 

The next two representative examples of special zeta functions serve 
to demonstrate the style of the new series representations we are able to 
obtain from the generalized generating function transformations established by 
this article. 
The \emph{Dirichlet beta function}, $\beta(s)$, is defined for 
$\Re(s) > 0$ by the series 
\begin{align*} 
\beta(s) & = \sum_{n \geq 0} \frac{(-1)^n}{(2n+1)^s} = 
     2^{-s} \Phi(-1, s, 1/2). 
\end{align*} 
The series in the previous equation is expanded through the generalized 
coefficients when $(\alpha, \beta) = (2, 1)$ as in the listings in 
Table \ref{table_TableOfGenCoeffs_ab21}. 
The first few special cases of $s$ over the positive integers are expanded 
by the following new series \citep[\cf \S 8]{FLAJOLETESUMS}: 
\begin{align*}
\beta(1) & = \sum_{j \geq 0} \frac{1}{2^{j+1}} 
     \binom{j+\frac{1}{2}}{\frac{1}{2}}^{-1} \\ 
\beta(2) & = \sum_{j \geq 0} \frac{1}{2^{j+1}} 
     \binom{j+\frac{1}{2}}{\frac{1}{2}}^{-1} \cdot \left( 
     1 + \HZNumber{j}{1}{2}{1}\right) \\ 
\beta(3) & = \sum_{j \geq 0} \frac{1}{2^{j+1}} 
     \binom{j+\frac{1}{2}}{\frac{1}{2}}^{-1} \cdot \left( 
     1 + \HZNumber{j}{1}{2}{1} + \frac{1}{2}\left( 
     \HZNumber{j}{1}{2}{1}^2 + \HZNumber{j}{2}{2}{1}\right) 
     \right). 
\end{align*} 
For $|z| < 1$, the \emph{Legendre chi function}, $\chi_\nu(z)$, 
is defined by the series 
\begin{align*} 
\chi_{\nu}(z) & = \sum_{k \geq 0} \frac{z^{2k+1}}{(2k+1)^\nu} = 
     2^{-\nu} z \times \Phi(z^2, \nu, 1/2). 
\end{align*} 
The first few positive integer cases of $\nu \geq 1$ are similarly 
expanded by the forms of the next series given by 
\begin{align*}
\chi_1(z) & = \sum_{j \geq 0} 
     \binom{j+\frac{1}{2}}{\frac{1}{2}}^{-1} 
     \frac{z \cdot (-z^2)^j}{(1-z^2)^{j+1}} \\ 
\chi_2(z) & = \sum_{j \geq 0} 
     \binom{j+\frac{1}{2}}{\frac{1}{2}}^{-1} \cdot \left( 
     1 + \HZNumber{j}{1}{2}{1}\right) \frac{z \cdot (-z^2)^j}{(1-z^2)^{j+1}} \\ 
\chi_3(z) & = \sum_{j \geq 0} 
     \binom{j+\frac{1}{2}}{\frac{1}{2}}^{-1} \cdot \left( 
     1 + \HZNumber{j}{1}{2}{1} + \frac{1}{2}\left( 
     \HZNumber{j}{1}{2}{1}^2 + \HZNumber{j}{2}{2}{1}\right) 
     \right) \frac{z \cdot (-z^2)^j}{(1-z^2)^{j+1}}. 
\end{align*} 
Other special function series and 
motivating examples of these generalized generating function transformations 
we consider within the examples in Section \ref{Section_Examples} of the 
article include special cases of the 
Hurwitz zeta and Lerch transcendent function series in 
\eqref{eqn_GenModifiedHZetaFn_def} and \eqref{eqn_ModLTransFn_series_def}. 
In particular, we consider concrete new series expansions of 
BBP-like series for special constants, 
\emph{alternating} and \emph{exotic Euler sums} with cubic denominators, 
\emph{alternating zeta function} sums with quadratic denominators, 
\emph{polygamma functions}, and 
several particular series defining the \emph{Riemann zeta function}, 
$\zeta(2k+1)$, over the odd positive integers.

\section{Generalized Stirling Numbers of the First Kind} 
\label{subSection_GenHZStirlingNumbers} 

\begin{table}[ht] 
\renewcommand{\arraystretch}{1} 
\setlength{\arraycolsep}{3pt} 
\begin{equation*} 
\begin{array}{|c|lllllllll|} \hline 
\trianglenk{j}{k} & 0 & 1 & 2 & 3 & 4 & 5 & 6 & 7 & 8 \\ \hline 
0 & 1 & 0 & 0 & 0 & 0 & 0 & 0 & 0 & 0 \\
1 & 0 & 1 & 0 & 0 & 0 & 0 & 0 & 0 & 0 \\
2 & 0 & 3 & 1 & 0 & 0 & 0 & 0 & 0 & 0 \\
3 & 0 & 15 & 8 & 1 & 0 & 0 & 0 & 0 & 0 \\
4 & 0 & 105 & 71 & 15 & 1 & 0 & 0 & 0 & 0 \\
5 & 0 & 945 & 744 & 206 & 24 & 1 & 0 & 0 & 0 \\
6 & 0 & 10395 & 9129 & 3010 & 470 & 35 & 1 & 0 & 0 \\
7 & 0 & 135135 & 129072 & 48259 & 9120 & 925 & 48 & 1 & 0 \\
8 & 0 & 2027025 & 2071215 & 852957 & 185059 & 22995 & 1645 & 63 & 1 \\ \hline 
\end{array}
\end{equation*} 
\caption{The Generalized Stirling Numbers of the First Kind, $\HZCf{k}{j}{2}{1}$} 
\label{table_TableOfGenS1Coeffs_kj_ab21}
\end{table} 

\subsection{Definition and Generating Functions} 

We first define a generalized set of coefficients in the symbolic 
polynomial expansions of the next products over $x$ as an extension of the 
results first given in \citep{MULTIFACTJIS}\footnote{ 
     The notation for \textit{Iverson's convention}, 
     $\Iverson{n = k} = \delta_{n,k}$, is consistent with its usage in 
     \citep{GKP}. 
}. 
\begin{align}
\label{eqn_HZCf_finite_product_def} 
\HZCf{n}{k}{\alpha}{\beta} & = 
    [x^{k-1}] x(x+\alpha+\beta)(x+2\alpha+\beta) \cdots 
    (x + (n-1) \alpha + \beta) \Iverson{n \geq 1} 
\end{align}
The polynomial coefficients of the powers of $x$ in 
\eqref{eqn_HZCf_finite_product_def}
are then defined by the following triangular recurrence for 
natural numbers $n, k \geq 0$: 
\begin{align} 
\label{eqn_gen_triangle_rdef_v1}
\HZCf{n}{k}{\alpha}{\beta} & = 
     (\alpha n+\beta-\alpha) \HZCf{n-1}{k}{\alpha}{\beta} + 
     \HZCf{n-1}{k-1}{\alpha}{\beta} + \Iverson{n = k = 0}. 
\end{align}
We also easily arrive at generating functions for the 
column sequences and for the generalized analogs to the 
\emph{Stirling convolution polynomials}, $\sigma_n(x)$ and 
$\sigma_n^{(\alpha)}(x)$, defined by 
\begin{align} 
\notag 
\sigma_n^{(\alpha, \beta)}(x) & = \HZCf{x}{x-n}{\alpha}{\beta} 
     \frac{(x-n-1)!}{x!}. 
\end{align} 
The series enumerating these coefficients are expanded as the 
following closed--form generating functions 
\citep[\S 6.2]{GKP} \citep{CVLPOLYS}\footnote{ 
     For fixed $x, \alpha, \beta$, we have a known identity for the 
     following exponential generating functions, which then implies the 
     first result in \eqref{eqn_HZCfLogRowGFs_GenSPoly_GFs} by 
     considering powers of $x^k$ as functions of $z$ where 
     $(x)_n$ denotes the \emph{Pochhammer symbol} 
     \citep{CVLPOLYS} \citep[\S 5.2(iii)]{NISTHB}: 
     \begin{align*} 
     \sum_{n \geq 0} \left(\frac{x+\beta}{\alpha}\right)_n 
     \frac{(\alpha z)^n}{n!} & = (1-\alpha z)^{-(x+\beta) / \alpha}. 
     \end{align*} 
     We can then apply a double integral transformation for the 
     \emph{beta function}, $B(a, b)$, in the form of \citep[\S 5.12]{NISTHB} 
     \begin{align*} 
     B(a, b)^{2} & = \binom{a+b}{a}^{-2} \cdot \frac{(a+b)^2}{a^2 b^2} = 
          \int_0^{\infty} \int_0^{\infty} 
          \frac{(ts)^{a-1}}{\left[(1+t)(1+s)\right]^{a+b}} dt\ ds, 
     \end{align*} 
     for real numbers $a, b > 0$ such that 
     $b \in \mathbb{Q} \setminus \mathbb{Z}$ to these generating functions 
     to obtain partially complete integral representations for the 
     transformed series over the modified coefficients in 
     \eqref{eqn_HNumExps_OfThe_HZetaCoeffs} and 
     \eqref{eqn_S2SStarCfhab_HNum_exp_rec_formula_v2} -- 
     \eqref{eqn_S2SStarCfhab_HNum_exp_formula_finite_sum_v1} of 
     Section \ref{subSection_Intro_GenGFSeriesEnums} cited by the examples in 
     Section \ref{Section_Examples}. 
}: 
\begin{align} 
\label{eqn_HZCfLogRowGFs_GenSPoly_GFs} 
\sum_{n \geq 0} \HZCf{n}{k}{\alpha}{\beta} \frac{z^n}{n!} & = 
     \frac{(1-\alpha z)^{-\beta / \alpha}}{k! \alpha^k} 
     \Log(1-\alpha z)^k \\ 
\notag 
\sum_{n \geq 0} x \sigma_n^{(\alpha, \beta)}(x) z^n & = 
     e^{\beta z} \left(\frac{\alpha z e^{\alpha z}}{e^{\alpha z} - 1} 
     \right)^x. 
\end{align} 
When $(\alpha, \beta) = (1, 0), (\alpha, 1-\alpha)$ we arrive at the 
definitions of the respective triangular recurrences defining the 
\emph{Stirling numbers of the first kind}, $\gkpSI{n}{k}$, and the generalized 
\emph{$\alpha$--factorial functions}, $n!_{(\alpha)}$, from the 
references \citep{STIRESUMS,GKP,MULTIFACTJIS}. 
Table \ref{table_TableOfGenS1Coeffs_kj_ab21} 
lists the first several rows of the triangle in \eqref{eqn_gen_triangle_rdef_v1}
corresponding to the special case of $(\alpha, \beta) = (2, 1)$ 
as considered in many special case expansions. 

\subsection{Expansions by the Generalized Harmonic Number Sequences} 

We find, as in the references \citep{STIRESUMS,GFTRANS2016}, 
that this generalized form of a Stirling--number--like 
triangle satisfies a number of analogous 
harmonic number expansions to the Stirling numbers of the first kind 
given in terms of the partial sums, 
$\HZNumber{n}{r}{\alpha}{\beta} = \sum_{k=1}^n (\alpha k+\beta)^{-r}$, 
of the \emph{modified Hurwitz zeta function}, 
$\zeta(s, \alpha, \beta) = \sum_{n \geq 1} 1 / (\alpha n + \beta)^s$. 
For example, we may expand special case formulas for the 
triangle columns at $k = 2, 3, 4$ in the following forms\footnote{
     We define the shorthand factorial function notation 
     for the products as 
     $\HZFact{n}{\alpha}{\beta} = \prod_{j=1}^n (\alpha j+\beta)$. 
}: 
\begin{align}
\label{eqn_}
\HZCf{n+1}{2}{\alpha}{\beta} & = \HZFact{n}{\alpha}{\beta} \times 
     \HZNumber{n}{1}{\alpha}{\beta} \\ 
\notag 
\HZCf{n+1}{3}{\alpha}{\beta} & = \frac{\HZFact{n}{\alpha}{\beta}}{2} \times 
     \left(\left(\HZNumber{n}{1}{\alpha}{\beta}\right)^2 - 
     \HZNumber{n}{2}{\alpha}{\beta}\right) \\ 
\notag 
\HZCf{n+1}{4}{\alpha}{\beta} & = \frac{\HZFact{n}{\alpha}{\beta}}{6} \times 
     \left(\left(\HZNumber{n}{1}{\alpha}{\beta}\right)^3 - 
     3 \HZNumber{n}{1}{\alpha}{\beta} \HZNumber{n}{2}{\alpha}{\beta} + 
     2 \HZNumber{n}{3}{\alpha}{\beta}\right). 
\end{align}
Similarly, we invert to expand the first few cases of the 
generalized $r$--order harmonic numbers through products of the coefficients in 
\eqref{eqn_gen_triangle_rdef_v1} as \citep[\cf \S 4.3]{MULTIFACTJIS} 
\begin{align} 
\notag 
\HZNumber{n}{2}{\alpha}{\beta} & = 
     \frac{1}{\left(\HZFact{n}{\alpha}{\beta}\right)^2}\left(
     \HZCf{n+1}{2}{\alpha}{\beta}^2 - 
     2 \HZCf{n+1}{1}{\alpha}{\beta} \HZCf{n+1}{3}{\alpha}{\beta}
     \right) \\ 
\notag 
\HZNumber{n}{3}{\alpha}{\beta} & = 
     \frac{1}{\left(\HZFact{n}{\alpha}{\beta}\right)^3}\Biggl(
     \HZCf{n+1}{2}{\alpha}{\beta}^3 - 
     3 \HZCf{n+1}{1}{\alpha}{\beta} \HZCf{n+1}{2}{\alpha}{\beta} 
     \HZCf{n+1}{3}{\alpha}{\beta} \\ 
\notag 
     & \phantom{= \frac{1}{\left(\HZFact{n}{\alpha}{\beta}\right)^3}\Biggl(
     \HZCf{n+1}{2}{\alpha}{\beta}^3\ } + 
     3 \HZCf{n+1}{1}{\alpha}{\beta}^2 \HZCf{n+1}{4}{\alpha}{\beta}
     \Biggr) \\ 
\notag 
\HZNumber{n}{4}{\alpha}{\beta} & = 
     \frac{1}{\left(\HZFact{n}{\alpha}{\beta}\right)^4}\Biggl(
     \HZCf{n+1}{2}{\alpha}{\beta}^4 - 
     4 \HZCf{n+1}{1}{\alpha}{\beta} \HZCf{n+1}{2}{\alpha}{\beta}^2 
     \HZCf{n+1}{3}{\alpha}{\beta} \\ 
\notag 
     & \phantom{= \frac{1}{\left(\HZFact{n}{\alpha}{\beta}\right)^4} 
                \Biggl(\qquad\ } + 
     2 \HZCf{n+1}{1}{\alpha}{\beta}^2 \HZCf{n+1}{3}{\alpha}{\beta}^2 - 
     4 \HZCf{n+1}{1}{\alpha}{\beta}^3 \HZCf{n+1}{5}{\alpha}{\beta} \\ 
\notag 
     & \phantom{= \frac{1}{\left(\HZFact{n}{\alpha}{\beta}\right)^4} 
                \Biggl(\qquad\ } + 
     4 \HZCf{n+1}{1}{\alpha}{\beta}^2 
     \HZCf{n+1}{2}{\alpha}{\beta} \HZCf{n+1}{4}{\alpha}{\beta} 
     \Biggr). 
\end{align} 
In general, we can use the \emph{elementary symmetric polynomials} 
implicit to the product-based definition of these generalized 
Stirling numbers in \eqref{eqn_HZCf_finite_product_def} to show that 
\begin{align*} 
\HZCf{n+1}{k}{\alpha}{\beta} & = 
     (-1)^k \cdot \HZFact{n}{\alpha}{\beta} \times 
     Y_k\left(-\HZNumber{n}{1}{\alpha}{\beta}, 
     \ldots, (-1)^k \HZNumber{n}{k}{\alpha}{\beta} \cdot (k-1)!\right) \\ 
\HZNumber{n}{k}{\alpha}{\beta} & = (-1)^k (k+1) 
     \HZCf{n+1}{1}{\alpha}{\beta}^{k+1} [t^{k+1}] 
     \Log\left(\HZCf{n+1}{1}{\alpha}{\beta} + 
     \sum_{j \geq 1} \HZCf{n+1}{j+1}{\alpha}{\beta} t^j 
     \right), 
\end{align*} 
where $Y_n(x_1, x_2, \ldots, x_n)$ denotes the 
\emph{exponential}, or \emph{complete}, \emph{Bell polynomial} whose 
exponential generating function is given by 
$\Phi(t, 1) \equiv \exp\left(\sum_{j \geq 1} x_j t^j / j!\right)$ 
\citep[\S 4.1.8]{UC}. 

Additionally, the next recurrences are obtained for the 
generalized harmonic numbers in terms of these coefficients 
corresponding to the partial sums in our definition of the 
``\emph{modified}'' Hurwitz zeta function, 
$\zeta(s, \alpha, \beta) = \alpha^{-s} \times \zeta(s, \beta / \alpha)$. 
\begin{align} 
\notag 
\HZNumber{n}{p}{\alpha}{\beta} & = 
     \sum_{1 \leq j < p} \HZCf{n+1}{p+1-j}{\alpha}{\beta} 
     \frac{(-1)^{p+1-j} \HZNumber{n}{j}{\alpha}{\beta}}{ 
     \HZFact{n}{\alpha}{\beta}} + 
     \HZCf{n+1}{p+1}{\alpha}{\beta} \frac{p (-1)^{p+1}}{ 
     \HZFact{n}{\alpha}{\beta}} \\ 
\notag 
\HZNumber{n+1}{p}{\alpha}{\beta} & = \HZNumber{n}{p}{\alpha}{\beta} + 
     \HZCf{n+1}{p}{\alpha}{\beta} \frac{(-1)^{p+1}}{
     \HZFact{(n+1)}{\alpha}{\beta}} \\ 
\notag 
     & \phantom{= \HZNumber{n}{p}{\alpha}{\beta}\ } + 
     \sum_{1 \leq j < p} \HZCf{n+2}{p+1-j}{\alpha}{\beta} 
     \frac{(-1)^{p+1-j}}{(\alpha n + \alpha + \beta)^j 
     \HZFact{(n+1)}{\alpha}{\beta}} \\ 
\notag 
\HZNumber{n+1}{p}{\alpha}{\beta} & = \HZNumber{n}{p}{\alpha}{\beta} + 
     \frac{1}{(\alpha n+\alpha+\beta)^{p-1}} \\ 
\notag 
     & \phantom{= \HZNumber{n}{p}{\alpha}{\beta}\ } + 
     \frac{(-1)^{p-1}}{\HZFact{(n+1)}{\alpha}{\beta}}\left(
     \HZCf{n+1}{p}{\alpha}{\beta} + \HZCf{n+1}{p-1}{\alpha}{\beta}\right) \\ 
\notag
     & \phantom{= \HZNumber{n}{p}{\alpha}{\beta}\ } + 
     \HZCf{n+2}{p}{\alpha}{\beta} \frac{(-1)^p}{(\alpha n+\alpha+\beta) 
     \HZFact{(n+1)}{\alpha}{\beta}} \\ 
\notag 
     & \phantom{= \HZNumber{n}{p}{\alpha}{\beta}\ } + 
     \sum_{j=0}^{p-3} \HZCf{n+2}{j+2}{\alpha}{\beta} 
     \frac{(-1)^{j+1} (\alpha n+\alpha + \beta - 1)}{ 
     (\alpha n + \alpha + \beta)^{p-1-j} \HZFact{(n+1)}{\alpha}{\beta}} 
\end{align} 
The last equation provides an implicit functional equation between our 
modified Hurwitz zeta function involving the generalized 
Stirling numbers of the first kind in \eqref{eqn_gen_triangle_rdef_v1}. 
For $(\alpha, \beta) := (1, 0)$, the previous equation implies 
new functional equations relating the $p$-order and $(p-1)$-order 
polylogarithm functions, $\Li_s(z) \equiv \Phi(z, s, 1, 0)$. 

\section{Transformations of Ordinary Power Series by 
         Generalized Stirling Numbers of the Second Kind} 
\label{subSection_Intro_GenGFSeriesEnums} 

\subsection{Definitions and Preliminary Examples} 
\label{subSection_GenCoeffs_SpCases_ab_21_31_32} 

Another approach to the relations of the generalized harmonic number 
sequences to the forms of the triangles defined by 
\eqref{eqn_gen_triangle_rdef_v1} proceeds as in the reference 
\citep{GFTRANS2016}. 
In particular, the next definitions lead to new expansions of 
many series and BBP--type formulas for special functions and constants 
from the introduction and in Section \ref{Section_Examples} which are 
implied by the new identities we prove for the series expansions of the 
\emph{modified Lerch transcendent} function, 
$\Phi(z, s, \alpha, \beta) = \alpha^{-s} \Phi(z, s, \beta / \alpha)$. 
\begin{align} 
\label{eqn_HZCfStar_sum_def_v1} 
\HZCfStar{k}{j}{\alpha}{\beta} & = \frac{1}{j!} \sum_{0 \leq m \leq j} 
     \binom{j}{m} \frac{(-1)^{j-m}}{(\alpha m+\beta)^{k-2}} \\ 
\notag 
\Phi(z, s, \alpha, \beta) & = \beta^{-s} + 
     \sum_{j \geq 0} \HZCfStar{s+2}{j}{\alpha}{\beta} \frac{j! z^j}{(1-z)^{j+1}} 
\end{align} 
The definition of the generalized 
\emph{Stirling numbers of the second kind}\footnote{ 
     See the conclusions in Section \ref{subSection_Concl_S2Analog} 
     for a short discussion of why we consider these transformation 
     coefficients to be generalized Stirling numbers of the second kind. 
} 
provided by \eqref{eqn_HZCfStar_sum_def_v1}, and recursively by 
\begin{align*} 
\HZCfStar{k}{j}{\alpha}{\beta} & = 
     (\alpha j + \beta) \HZCfStar{k+1}{j}{\alpha}{\beta} + 
     \alpha \cdot \HZCfStar{k+1}{j-1}{\alpha}{\beta}, 
\end{align*} 
also implies the next new truncated partial power series identities 
for the ``\emph{modified}'' Lerch transcendent function over 
some sequence, $\langle g_n \rangle$, whose ordinary generating function, 
$G(z)$, has derivatives of all orders, 
for $\alpha \geq 1$ and $0 \leq \beta < \alpha$, and 
for any fixed $u \geq 1, u_0 \geq 0$ (See \citep{GFTRANS2016}). 
\begin{align} 
\label{eqn_} 
\sum_{n=1}^{u} \frac{g_n}{(\alpha n + \beta)^{k}} z^n & = [w^u]\left( 
     \sum_{j=1}^{u+u_0} \HZCfStar{k+2}{j}{\alpha}{\beta}  
     \frac{(wz)^{j} G^{(j)}(wz)}{(1-w)} \right) \\ 
\notag 
\sum_{1 \leq n \leq u} \HZNumber{n}{k}{\alpha}{\beta} z^n & = [w^u]\left( 
     \sum_{1 \leq j \leq u+u_0} \HZCfStar{k+2}{j}{\alpha}{\beta} 
     \frac{(wz)^{j} \cdot j!}{(1-w) (1-wz)^{j+2}} \right) \\ 
\notag 
\sum_{1 \leq n \leq u} \HZNumber{n}{k}{\alpha}{\beta} \frac{z^n}{n!} & = 
     [w^u]\left( 
     \sum_{1 \leq j \leq u+u_0} \HZCfStar{k+2}{j}{\alpha}{\beta} 
     \frac{(wz)^{j} \cdot e^{wz}\left(j+1 + wz\right)}{(j+1) (1-w)} \right) \\ 
\notag 
\sum_{1 \leq n \leq u} \left(\sum_{k=1}^{n} \frac{t^k}{(\alpha k + \beta)^{r}} 
     \right) z^n & = [w^u]\left( 
     \sum_{1 \leq j \leq u} \HZCfStar{r+2}{j}{\alpha}{\beta}  
     \frac{(twz)^{j} \cdot j!}{(1-w) (1-wz) (1-twz)^{j+1}}\right)
\end{align} 
The generalized harmonic number expansions of the coefficients in 
\eqref{eqn_HZCfStar_sum_def_v1} are considered next in 
Section \ref{subSection_HNumExps_OfThe_HZetaCoeffs}. 
An even more general proof of the 
formal power series transformation suggested in the 
concluding remarks from \citep{GFTRANS2016} is given below in 
Section \ref{Section_CombProofOfGenfnPowk_Formal_SeriesTrans}. 
Table \ref{table_TableOfGenCoeffs_ab21}, 
Table \ref{table_TableOfGenCoeffs_ab31}, and 
Table \ref{table_TableOfGenCoeffs_ab32} 
each provide listings of useful particular special cases of the 
generalized transformation coefficients, or alternately, 
\emph{generalized Stirling numbers of the second kind} 
within the context of this article, corresponding to 
$(\alpha, \beta) := (2, 1), (3, 1), (3, 2)$, respectively. 

\begin{table}[t] 
\renewcommand{\arraystretch}{1.25} 
\setlength{\arraycolsep}{3pt} 
\begin{equation*} 
\begin{array}{|c|lllllll|} \hline 
\trianglenk{j}{k} & 0 & 1 & 2 & 3 & 4 & 5 & 6 \\ \hline 
0 & 0 & 0 & 0 & 0 & 0 & 0 & 0 \\
1 & 9 & 3 & 1 & \frac{1}{3} & \frac{1}{9} & \frac{1}{27} & \frac{1}{81} \\
2 & 7 & 1 & 1 & \frac{7}{15} & \frac{41}{225} & \frac{223}{3375} &
        \frac{1169}{50625} \\
3 & 1 & 1 & 1 & \frac{19}{35} & \frac{859}{3675} & 
    \frac{34739}{385875} & \frac{1323019}{40516875} \\
4 & 1 & 1 & 1 & \frac{187}{315} & \frac{27161}{99225} &
    \frac{3451843}{31255875} & \frac{406586609}{9845600625} \\
5 & 1 & 1 & 1 & \frac{437}{693} & \frac{735197}{2401245} &
    \frac{1066933061}{8320313925} & \frac{1418417467373}{28829887750125} \\
6 & 1 & 1 & 1 & \frac{1979}{3003} & \frac{45087479}{135270135} &
    \frac{877474863971}{6093243231075} &
    \frac{15505503106933439}{274470141343773375} \\
7 & 1 & 1 & 1 & \frac{4387}{6435} & \frac{103349119}{289864575} &
    \frac{2065307132299}{13056949780875} &
    \frac{1488524941286431}{23526012115180575} \\
8 & 1 & 1 & 1 & \frac{76627}{109395} & \frac{31562623583}{83770862175} &
    \frac{10971718559046811}{64148794273438875} &
    \frac{683894055421671560539}{9824580289359984022875} \\
    \hline 
\end{array}
\end{equation*} 
\caption{A Table of the Generalized Coefficients 
         $\HZCfStar{k}{j}{2}{1} \times j! (-1)^{j-1}$} 
\label{table_TableOfGenCoeffs_ab21}
\end{table} 

\begin{table}[ht] 
\renewcommand{\arraystretch}{1.25} 
\setlength{\arraycolsep}{3pt} 
\begin{equation*} 
\begin{array}{|c|lllllll|} \hline 
\trianglenk{j}{k} & 0 & 1 & 2 & 3 & 4 & 5 & 6 \\ \hline 
0 & 0 & 0 & 0 & 0 & 0 & 0 & 0 \\
1 & 16 & 4 & 1 & \frac{1}{4} & \frac{1}{16} & \frac{1}{64} & \frac{1}{256} \\
2 & -17 & 1 & 1 & \frac{5}{14} & \frac{41}{392} & \frac{311}{10976} &
    \frac{2273}{307328} \\
3 & 1 & 1 & 1 & \frac{59}{140} & \frac{2671}{19600} & \frac{107369}{2744000} &
    \frac{4060291}{384160000} \\
4 & 1 & 1 & 1 & \frac{212}{455} & \frac{133849}{828100} &
    \frac{73174943}{1507142000} & \frac{37005870001}{2742998440000} \\
5 & 1 & 1 & 1 & \frac{727}{1456} & \frac{1936973}{10599680} &
    \frac{4393719979}{77165670400} & \frac{9104269630637}{561766080512000} \\
6 & 1 & 1 & 1 & \frac{7271}{13832} & \frac{384155263}{1913242240} &
    \frac{17071846526411}{264639666636800} &
    \frac{686298711281124727}{36604958689202176000} \\
7 & 1 & 1 & 1 & \frac{23789}{43472} & \frac{14322370919}{66143517440} &
    \frac{7187615461845233}{100638684655308800} &
    \frac{3237486239486747349191}{153123771476745445376000} \\
8 & 1 & 1 & 1 & \frac{76801}{135850} & \frac{238206415289}{1033492460000} &
    \frac{611558324636496331}{7862397238696000000} &
    \frac{1400156984227714635455249}{59813973233103689600000000} \\
    \hline 
\end{array}
\end{equation*} 
\caption{A Table of the Generalized Coefficients 
         $\HZCfStar{k}{j}{3}{1} \times j! (-1)^{j-1}$} 
\label{table_TableOfGenCoeffs_ab31}
\end{table} 

\begin{table}[t] 
\renewcommand{\arraystretch}{1.25} 
\setlength{\arraycolsep}{3pt} 
\begin{equation*} 
\begin{array}{|c|lllllll|} \hline 
\trianglenk{j}{k} & 0 & 1 & 2 & 3 & 4 & 5 & 6 \\ \hline 
0 & 0 & 0 & 0 & 0 & 0 & 0 & 0 \\
1 & 25 & 5 & 1 & \frac{1}{5} & \frac{1}{25} & \frac{1}{125} & \frac{1}{625} \\
2 & -14 & 2 & 1 & \frac{11}{40} & \frac{103}{1600} & \frac{899}{64000} &
    \frac{7567}{2560000} \\
3 & 4 & 2 & 1 & \frac{139}{440} & \frac{15757}{193600} &
    \frac{1609291}{85184000} & \frac{155016733}{37480960000} \\
4 & 4 & 2 & 1 & \frac{527}{1540} & \frac{446837}{4743200} &
    \frac{334869917}{14609056000} & \frac{233183599997}{44995892480000} \\
5 & 4 & 2 & 1 & \frac{1889}{5236} & \frac{28606807}{274156960} &
    \frac{378441183599}{14354858425600} &
    \frac{4602491925840703}{751620387164416000} \\
6 & 4 & 2 & 1 & \frac{19619}{52360} & \frac{61764761}{548313920} &
    \frac{4214471373881}{143548584256000} &
    \frac{10491182677877357}{1503240774328832000} \\
7 & 4 & 2 & 1 & \frac{66337}{172040} & \frac{4956449573}{41436866240} &
    \frac{7985964568560547}{249507946377536000} &
    \frac{466567679887167456041}{60095485932707810816000} \\
8 & 4 & 2 & 1 & \frac{110258}{279565} & \frac{43971566839}{350141519728} &
    \frac{4710810017671083829}{137042239547861648000} &
    \frac{3642461006944413986125043}{429096793431016946178944000} \\
    \hline 
\end{array}
\end{equation*} 
\caption{A Table of the Generalized Coefficients 
         $\HZCfStar{k}{j}{3}{2} \times j! (-1)^{j-1}$} 
\label{table_TableOfGenCoeffs_ab32}
\end{table} 

\subsection{Harmonic Number Expansions of the 
            Generalized Transformation Coefficients} 
\label{subSection_HNumExps_OfThe_HZetaCoeffs} 

For integers $\alpha \geq 1$ and $0 \leq \beta < \alpha$, 
let the first few special cases of the 
functions, $\widetilde{R}_k(\alpha, \beta; j)$, defined by 
\begin{align*} 
\widetilde{R}_k(\alpha, \beta; j) & := 
     \binom{j+\frac{\beta}{\alpha}}{j} 
     \sum_{m=1}^{j} \binom{j}{m} 
     \frac{(-1)^{m+1} \cdot \alpha m}{(\alpha m+\beta)^{k}} \Iverson{k \geq 2} + 
     \Iverson{k = 1} + \Iverson{k = 0}, 
\end{align*} 
be expanded in explicit formulas by the next equations. 
\begin{align} 
\label{eqn_HNumExps_OfThe_HZetaCoeffs} 
\widetilde{R}_1(\alpha, \beta; j) & = 1 \\ 
\notag 
\widetilde{R}_2(\alpha, \beta; j) & = H_j^{(1)}(\alpha, \beta)  \\ 
\notag 
\widetilde{R}_3(\alpha, \beta; j) & = \frac{1}{2}\left( 
     H_j^{(1)}(\alpha, \beta)^{2} + H_j^{(2)}(\alpha, \beta) 
     \right) \\ 
\notag 
\widetilde{R}_4(\alpha, \beta; j) & = \frac{1}{6}\left( 
     H_j^{(1)}(\alpha, \beta)^{3} + 3 H_j^{(1)}(\alpha, \beta) 
     H_j^{(2)}(\alpha, \beta) + 2 H_j^{(3)}(\alpha, \beta)
     \right) \\ 
\notag
\widetilde{R}_5(\alpha, \beta; j) & = \frac{1}{24}\Bigl( 
     H_j^{(1)}(\alpha, \beta)^{4} + 6 H_j^{(1)}(\alpha, \beta)^{2} 
     H_j^{(2)}(\alpha, \beta) + 3 H_j^{(2)}(\alpha, \beta)^{2} \\ 
\notag 
     & \phantom{= \frac{1}{24}\Bigl(\ } + 
     8 H_j^{(1)}(\alpha, \beta) H_j^{(3)}(\alpha, \beta) + 
     6 H_j^{(4)}(\alpha, \beta)
     \Bigr). 
\end{align}
For larger cases of $k > 5$, we employ the following heuristic to 
generate the harmonic number expansions of these functions, which for 
concrete special cases are easily obtained from \emph{Mathematica's} 
\texttt{Sigma} package: 
\begin{equation} 
\label{eqn_S2SCfhab_S1-like_HNum_exp_recursive_ident_stmt_v1} 
\widetilde{R}_{m}(\alpha, \beta; j) = \sum_{i=0}^{m-2} 
     \frac{\widetilde{R}_{m-2-i}(\alpha, \beta; j)}{(m-1)} \cdot 
     H_j^{(i+1)}(\alpha, \beta) + \Iverson{m = 1}. 
\end{equation}
We then define a vague analog to the harmonic number expansions in 
\citep{GFTRANS2016} through the expansions of these functions as 
\begin{align} 
\label{eqn_S2SStarCfhab_HNum_exp_rec_formula_v2} 
\HZCfStar{k+2}{j}{\alpha}{\beta} & = \HZCfStar{k+1}{j}{\alpha}{\beta} 
     \cdot \frac{1}{\beta} + 
     \binom{j+\beta / \alpha}{\beta / \alpha}^{-1} 
     \frac{(-1)^{j}}{\beta \cdot j!} \cdot 
     \widetilde{R}_{k}(\alpha, \beta; j) \\ 
\label{eqn_S2SStarCfhab_HNum_exp_formula_finite_sum_v1} 
     & = \frac{(-1)^{j-1}}{\beta^{k} \cdot j!} + 
     \sum_{m=0}^{k-1} 
     \binom{j+\beta / \alpha}{\beta / \alpha}^{-1} 
     \frac{(-1)^{j}}{j!} \cdot 
     \frac{\widetilde{R}_{m+1}(\alpha, \beta; j)}{\beta^{k-m}}. 
\end{align} 
Notice that the heuristic we used to generate more involved cases of the 
expansions for the functions, $\widetilde{R}_k(\alpha, \beta; j)$, 
imply recurrences for the $k$--order generalized harmonic number 
sequences in the following forms when $k \geq 3$: 
\begin{align} 
\notag 
H_n^{(k)}(\alpha, \beta) & = 
     \frac{H_n^{(k-1)}(\alpha, \beta)}{\beta} + \sum_{j=0}^{n} 
     \binom{n+1}{j+1} \binom{j+\beta/\alpha}{\beta/\alpha}^{-1} \left( 
     \frac{(-1)^{j}}{\beta} \cdot \widetilde{R}_{k}(\alpha, \beta; j) 
     \right) \\ 
\notag 
H_n^{(k)}(\alpha, \beta) & = 
     \frac{H_n^{(k-2)}(\alpha, \beta)}{\beta^2} \\ 
\notag 
     & \phantom{=\quad\ } + \sum_{j=0}^{n} 
     \binom{n+1}{j+1} \binom{j+\beta/\alpha}{\beta/\alpha}^{-1} (-1)^{j} 
     \cdot \left( 
     \frac{\widetilde{R}_{k}(\alpha, \beta; j)}{\beta} + 
     \frac{\widetilde{R}_{k-1}(\alpha, \beta; j)}{\beta^2} 
     \right). 
\end{align}
A pair of harmonic and Hurwitz zeta function related identities that follow 
from the generalized coefficient definitions in 
\eqref{eqn_HZCfStar_sum_def_v1} are obtained by 
similar methods from the reference \citep{GFTRANS2016} 
for $n \geq 1$ as follows: 
\begin{align*} 
\frac{1}{(\alpha n + \beta)^{k}} & = \sum_{0 \leq j \leq n} 
     \binom{n}{j} \HZCfStar{k+2}{j}{\alpha}{\beta} \cdot j! \\ 
     & = 
     \sum_{1 \leq m \leq n} \left( 
     \sum_{m \leq j \leq n} 
     \gkpSI{j}{m} \HZCfStar{k+2}{j}{\alpha}{\beta} (-1)^{j} 
     \right) (-1)^m n^m \\ 
H_n^{(k)}(\alpha, \beta) & = 
     \sum_{0 \leq j \leq n} \binom{n+1}{j+1} 
     \HZCfStar{k+2}{j}{\alpha}{\beta} \cdot j! \\ 
     & = 
     \sum_{0 \leq p \leq n+1} \left( 
     \sum_{0 \leq j \leq n} \gkpSI{j+1}{p} \HZCfStar{k+2}{j}{\alpha}{\beta} 
     \frac{(-1)^{j+1}}{(j+1)} 
     \right) (-1)^{p} \cdot (n+1)^{p}. 
\end{align*} 

\subsection{Proofs of the Zeta Series Transformations of 
            Formal Power Series} 
\label{Section_CombProofOfGenfnPowk_Formal_SeriesTrans}

\subsubsection{Proofs of the Geometric and Exponential Series Transformations} 

We claim that for any sequence, $\langle f(n) \rangle$, which is not 
identically zero for $n \geq 0$, 
any natural numbers $k \geq 1$, and a sequence, 
$\langle g_n \rangle$, whose ordinary generating function, $G(z)$, has 
derivatives of all orders, we have the (formal) series transformation 
\begin{align} 
\label{eqn_GenfZetaSeries_TransformSeries_stmt_v1} 
\sum_{n \geq 1} \frac{g_n}{f(n)^k} z^n & = 
     \sum_{j \geq 1} \SIIStarf{k+2}{j}{f} z^j G^{(j)}(z), 
\end{align} 
where the coefficients implicit to the right--hand--side series are 
defined by 
\begin{align} 
\label{eqn_S2StarfCf_finite_sum_def} 
\SIIStarf{k+2}{j}{f} & = \frac{1}{j!} \sum_{1 \leq m \leq j} 
     \binom{j}{m} \frac{(-1)^{j-m}}{f(m)^{k}}. 
\end{align} 
As in \citep{GFTRANS2016}, we primarily only work with these generalized 
series when the sequence generating function of $g_n$ is some 
variation of the geometric or exponential series. 
Therefore, for the content of our article, it suffices to prove the 
next two cases. 

\begin{proof}[Proof of the Geometric Series Case] 
Let $g_n \equiv 1$ so that its corresponding $j^{th}$ derivative is 
given by $G^{(j)}(z) = j! / (1-z)^{j+1}$. 
We proceed to expand the right--hand--side of 
\eqref{eqn_GenfZetaSeries_TransformSeries_stmt_v1} as follows: 
\begin{align*} 
\sum_{j \geq 1} \SIIStarf{k+2}{j}{f} \frac{j! z^j}{(1-z)^{j+1}} & = 
     \sum_{j \geq 1} \left( 
     \sum_{1 \leq m \leq j} \binom{j}{m} \frac{(-1)^m}{f(m)^k} 
     \right) \frac{(-z)^j}{(1-z)^{j+1}} \\ 
     & = 
     \sum_{m \geq 1} \left( 
     \sum_{j \geq m} \binom{j}{m} \frac{(-z)^j}{(1-z)^{j+1}} \right) 
     \frac{(-1)^m}{m! f(m)^k}. 
\end{align*} 
For a fixed $c \neq 1$, a known binomial sum identity gives that 
\begin{equation*} 
\sum_{j \geq m} \binom{j}{m} c^j = \frac{c^m}{(1-c)^{m+1}}, 
\end{equation*} 
which then implies that when $c \mapsto -z / (1-z)$ we have 
\begin{align*} 
\sum_{m \geq 1} \left( 
     \sum_{j \geq m} \binom{j}{m} \frac{(-z)^j}{(1-z)^{j+1}} \right) 
     \frac{(-1)^m}{m! f(m)^k} & = 
     \sum_{m \geq 1} (-z)^m \times \frac{(-1)^m}{m! f(m)^k}
     \qedhere
\end{align*} 
\end{proof} 

\begin{proof}[Proof of the Exponential Series Case] 
Let $g_n \equiv r^n / n!$ so that its corresponding $j^{th}$ derivative is 
given by $G^{(j)}(z) = r^j e^{rz}$ for all $j$. 
In this case, we proceed to expand the right--hand--side of 
\eqref{eqn_GenfZetaSeries_TransformSeries_stmt_v1} as 
\begin{align*} 
\sum_{j \geq 1} \SIIStarf{k+2}{j}{f} (r z)^j e^{rz} & = 
     \sum_{j \geq 1} \left( 
     \sum_{1 \leq m \leq j} \frac{(-1)^m}{m! (j-m)! f(m)^k} 
     \right) (-rz)^j e^{rz} \\ 
     & = 
     \sum_{m \geq 1} \left( 
     \sum_{j \geq m} \frac{(-rz)^j}{(j-m)!} \right) 
     \frac{(-1)^m}{m! f(m)^k} e^{rz} \\ 
     & = 
     \sum_{m \geq 1} \left(e^{-rz} (-r z)^m\right) 
     \frac{(-1)^m}{m! f(m)^k} e^{rz} \\ 
     & = 
     \sum_{m \geq 1} \frac{r^m z^m}{m! f(m)^k}. 
     \qedhere 
\end{align*} 
\end{proof} 

\subsubsection{Remarks on Symbolic Transformation Coefficient Identities} 

We observe that most of the identities formulated in 
\citep[\S 3]{GFTRANS2016} are easily restated as symbolic identities for the 
coefficients in \eqref{eqn_S2StarfCf_finite_sum_def}. 
If $f(m)$ is polynomial in $m$, by expanding $1/f(m)$ in partial fractions 
over linear factors of $m$, we arrive at sums over the coefficient forms in 
\eqref{eqn_HNumExps_OfThe_HZetaCoeffs} of 
Section \ref{subSection_HNumExps_OfThe_HZetaCoeffs} above. 
If $1/f(z)$ is a meromorphic function, we may alternately compute the 
symbolic coefficients in \eqref{eqn_S2StarfCf_finite_sum_def} by the next 
\textit{N\"{o}rlund--Rice integral} 
over a suitable contour given by \citep{FS-NORLUND-RICE-INTS}
\begin{align*}
\SIIStarf{k}{j}{f} \cdot j! & := 
\sum_{1 \leq m \leq j} \binom{j}{m} \frac{(-1)^{j-m}}{f(m)^k} \\ 
     & \phantom{:} = 
     \frac{j!}{2\pi\imath} \oint 
     \frac{f(z)^{-k}}{z(z-1)(z-2) \cdots (z-j)} dz. 
\end{align*}
We provide a brief overview of examples of several finite sum 
expansions which are also easily proved along the lines given in the 
reference, and then quickly move on to the particular cases of the 
series transformations at hand in this article. 
For example, if we let the \emph{$r$--order $f$--harmonic numbers}, 
$F_n^{(r)}(f) := \sum_{k=1}^n f(k)^{-r}$, be defined for a non-zero-valued 
function, $f(n)$, we may expand the coefficients in 
\eqref{eqn_S2StarfCf_finite_sum_def} as 
\begin{align*} 
\SIIStarf{k}{j}{f} & = 
     \sum_{0 \leq i < j} \frac{(j+1) (-1)^{j-1-i}}{(j-1-i)! (i+2)!} \times 
     F_{i+1}^{(k)}(f) \\ 
     & = 
     \sum_{0 \leq i < j} \frac{(j+1) (-1)^{j-1-i}}{(j-1-i)! (i+2)!} \left( 
     F_{i+2}^{(k)}(f) - \frac{1}{f(i+2)^k} \right) \\ 
     & = 
     \sum_{0 \leq i < j} \frac{(-1)^{j-1-i} F_{i+1}^{(k-r)}(f)}{ 
     (j-1-i)! (i+2)!} \left( 
     \frac{(i+2)}{f(i+1)^r} + \frac{(j-1-i)}{f(i+2)^r} \right) \\ 
\SIIStarf{k}{j}{f} & = 
     \sum_{m=0}^k \sum_{i=1}^k \sum_{r=1}^j 
     \gkpSI{k}{i} \gkpSII{i-1}{m} \binom{j}{r} 
     \frac{(-1)^{j-r+m} m! (f(r) - 1)^{m}}{j! (k-1)! f(r)^{m+1}}, 
\end{align*} 
for any integers $j, k \geq 1$ and real--valued $r \in (0, k)$. 
Notice that we can also similarly define the symbolic generalized 
Stirling numbers of the first kind by 
\begin{align*} 
\gkpSI{n}{k}_{f} & := [x^k] \prod_j \left(x+f(j)\right), 
\end{align*} 
and then proceed to derive a whole new related set of even more general 
symbolic combinatorial identities and properties for these coefficients 
involving the $f$-harmonic-numbers, $F_n^{(r)}(f)$. 

\section{Examples of New Series Expansions for Modified Zeta Function Series} 
\label{Section_Examples} 

\subsection{A Special Class of Alternating Euler Sums} 

A first pair of alternating Euler sums related to the 
Dirichlet beta function constants are expanded in the following forms 
\citep[\S 8]{FLAJOLETESUMS} \footnote{ 
     These special cases of the generalized harmonic number sequences are 
     expanded in terms of the ordinary $r$-order harmonic numbers, 
     $H_n^{(r)}$, considered in the expansions of \citep{GFTRANS2016} as 
     \begin{align*} 
     \HZNumber{n}{r}{2}{1} & = H_{2n+2}^{(r)} - 2^{-r} \cdot H_{n+1}^{(r)} - 1. 
     \end{align*} 
}: 
\begin{align*} 
\sum_{n \geq 0} \frac{(-1)^n}{(2n+1)} H_n & = 
     \sum_{n \geq 0} \frac{(-1)^n}{(2n+1)^2} - \frac{\pi}{2} \Log(2) \\ 
     & = 
     \sum_{j \geq 0} \binom{j+\frac{1}{2}}{\frac{1}{2}}^{-1} 
     \frac{\left(H_j-\Log(2)\right)}{2^{j+1}} \\ 
\sum_{n \geq 0} \frac{(-1)^n}{(2n+1)^3} H_n & = 
     3 \sum_{n \geq 0} \frac{(-1)^n}{(2n+1)^4} - \frac{7\pi}{16} \zeta(3) - 
     \frac{\pi^3}{16} \Log(2) \\ 
     & = 
     \sum_{j \geq 0} \binom{j+\frac{1}{2}}{\frac{1}{2}}^{-1} 
     \frac{\left(H_j-\Log(2)\right)}{2^{j+1}} \times \\ 
     & \phantom{= \sum\binom{j+1/2}{1/2}\ } \times 
     \left( 
     1 + \HZNumber{j}{1}{2}{1} + \frac{1}{2}\left(\HZNumber{j}{1}{2}{1}^2 + 
     \HZNumber{j}{2}{2}{1}\right)
     \right). 
\end{align*} 
The second and fourth series on the right--hand--side of the 
previous equations follow from an identity for the $j^{th}$ derivatives 
of the first--order harmonic number generating function given by 
\begin{align*} 
D_z^{(j)}\left[-\frac{\Log(1-z)}{(1-z)}\right] & = 
     \frac{\left(H_j-Log(1-z)\right) j!}{(1-z)^{j+1}}, 
\end{align*} 
for all integers $j \geq 0$. 
Since this identity is straightforward to 
prove by induction, we move quickly along to the next example. 

We also observe the key difference between the 
generalized zeta series transform coefficients introduced in 
\citep{GFTRANS2016} and those defined by 
\eqref{eqn_HNumExps_OfThe_HZetaCoeffs} and 
\eqref{eqn_S2SStarCfhab_HNum_exp_formula_finite_sum_v1} of this article. 
In particular, the definitions of the generalized coefficients given in this 
article imply functional equations for a number of special series. 
For example, the second series in the equations immediately above is given in 
terms of the first sum in the following form: 
\begin{align*} 
\sum_{n \geq 0} \frac{(-1)^n}{(2n+1)^3} H_n & = 
     \sum_{n \geq 0} \frac{(-1)^n}{(2n+1)^2} H_n \\ 
     & \phantom{=\ } + 
     \sum_{j \geq 0} \binom{j+\frac{1}{2}}{\frac{1}{2}}^{-1} 
     \frac{(H_j - \Log(2))}{2^{j+2}} \left( 
     \HZNumber{j}{1}{2}{1}^2 + \HZNumber{j}{2}{2}{1}\right). 
\end{align*} 
We have similar relations for series defining rational multiples of the 
\emph{polygamma functions}, $\psi^{s-1}(z/2) - \psi^{s-1}((z+1)/2)$, 
for example, as in the next pair of related sums given by 
\begin{align*} 
\sum_{k \geq 0} \frac{(-1)^k}{(k+z)^2} & = 
     \sum_{j \geq 0} \binom{j+z}{z}^{-1} \left[ 
     \frac{1}{z^2} + \frac{1}{z} \HZNumber{j}{1}{1}{z} 
     \right] \frac{1}{2^{j+1}} \\ 
\sum_{k \geq 0} \frac{(-1)^k}{(k+z)^3} & = 
     \frac{1}{z}\left(\sum_{k \geq 0} \frac{(-1)^k}{(k+z)^2}\right) + 
     \sum_{j \geq 0} \binom{j+z}{z}^{-1} \left[ 
     \frac{1}{z}\left(\HZNumber{j}{1}{1}{z}^2 + \HZNumber{j}{2}{1}{z}\right) 
     \right] \frac{1}{2^{j+2}}, 
\end{align*} 
which then implies functional equations between the polygamma functions 
such as the following identity: 
\begin{align*} 
-\frac{z}{4} & \left(\psi^{(2)}\left(\frac{z}{2}\right) - 
     \psi^{(2)}\left(\frac{z+1}{2}\right)\right) = \\ 
     & \phantom{\quad\ } 
     \psi^{\prime}\left(\frac{z}{2}\right) - 
     \psi^{\prime}\left(\frac{z+1}{2}\right) + 4 \times 
     \sum_{j \geq 0} \binom{j+z}{z}^{-1} \left[ 
     \HZNumber{j}{1}{1}{z}^2 + \HZNumber{j}{2}{1}{z} 
     \right] \frac{1}{2^{j+2}}. 
\end{align*} 

\subsection{An Exotic Euler Sum with Powers of Cubic Denominators} 

Flajolet mentions a more ``\emph{exotic}'' family of Euler sums in his 
article defined for positive integers $q$ by \citep{FLAJOLETESUMS} 
\begin{equation*} 
A_q^{\ast} = \sum_{n \geq 1} \frac{(-1)^n H_n^2}{\left[(2n-1)(2n)(2n+1)\right]^q}. 
\end{equation*} 
It is not difficult to prove that we have the following two 
ordinary generating functions for the squares of the 
first--order harmonic numbers: 
\begin{align} 
\label{eqn_HnSquared_GFs} 
\sum_{n \geq 0} H_n^2 z^n & = 
     \frac{1}{(1-z)}\left(\Log(1-z)^2 + \Li_2(z)\right) \\ 
\notag 
     & = 
     -\frac{1}{(1-z)}\left(2 \Li_2\left(-\frac{z}{1-z}\right) + \Li_2(z) 
     \right). 
\end{align} 
A proof of these two series identities follows from the expansions of the 
\emph{polylogarithm function}, $\Li_2(z) / (1-z)$, generating the 
second--order harmonic numbers, $H_n^{(2)}$, in \citep[\S 4]{GFTRANS2016}. 
In particular, we expand the polylogarithm series as 
\begin{align*}
\frac{\Li_2(z)}{(1-z)} & = - \sum_{j \geq 0} 
     \frac{\left(H_j^2+H_j^{(2)}\right)}{2 (1-z)^2} 
     \left(-\frac{z}{1-z}\right)^j, 
\end{align*} 
and then perform the change of variable $z \mapsto -z / (1-z)$ to 
obtain these results. 

Since the derivatives of the polylogarithm functions in each of the 
equations in \eqref{eqn_HnSquared_GFs} are tedious and messy to expand, 
we do not give any explicit series for these Euler sums, $A_q^{\ast}$. 
However, we do note that a sufficiently motivated reader may expand these 
sums by the generalized coefficients we defined in 
\eqref{eqn_HZCfStar_sum_def_v1} by taking partial fractions of the 
denominators of $A_q$. 
For example, when $q = 1, 2$ we have the series 
\begin{align*} 
A_1^{\ast} & = \sum_{n \geq 1} \frac{(-1)^n H_n^2}{2} \left(\frac{1}{(2n+1)} - 
     \frac{1}{n} + \frac{1}{(2n-1)}\right) \\ 
A_2^{\ast} & = \sum_{n \geq 1} \frac{(-1)^n H_n^2}{4} \left(\frac{3}{(2n+1)} - 
     \frac{3}{(2n-1)} + \frac{1}{(2n+1)^2} + 
     \frac{1}{n^2} + \frac{1}{(2n-1)^2}\right). 
\end{align*} 
To expand the more involved cases of the polylogarithm function derivatives, 
we first note that for integers $s \geq 2$ and $|z| \leq 1$ we have 
\citep[\S 2.7]{SPECIAL-FUNCTIONS}
\begin{equation*} 
D_z\left[\Li_s(z)\right] = \frac{1}{z} \Li_{s-1}(z), 
\end{equation*} 
where for $r \in \mathbb{Z}^{+}$ we have the identity that 
\citep[\S 7.4]{GKP} 
\begin{align*} 
\Li_{-r}(z) & = \sum_{0 \leq j \leq r} \gkpSII{r}{j} \frac{z^j j!}{(1-z)^{j+1}} = 
     \frac{1}{(1-z)^{r+1}} \times \sum_{0 \leq i \leq r} \gkpEI{r}{i} z^{i+1}, 
\end{align*} 
and where the composite derivatives in the second equation of 
\eqref{eqn_HnSquared_GFs} are expanded by \emph{Fa\'{a} de Bruno's formula} 
\citep[\S 1.4(iii)]{NISTHB}. 

\subsection{Zeta Function Series with Powers of Quadratic Denominators}

We next provide examples of generalized zeta function series over 
denominator powers of quadratic polynomials. 
The method employed to expand these particular series is to factor and then 
take partial fractions to apply the generalized transformation cases 
we study within this article. 
For example, we observe that 
\begin{align*} 
\frac{1}{n^2+1} & = \frac{\imath}{2}\left( 
     \frac{1}{n+\imath} - \frac{1}{n-\imath}\right) \\ 
\frac{1}{(n^2+1)^2} & = -\frac{1}{4}\left(\frac{\imath}{(n-\imath)} - 
     \frac{\imath}{(n+\imath)} + \frac{1}{(n+\imath)^2} + 
     \frac{1}{(n-\imath)^2} \right), 
\end{align*} 
which immediately leads to the first two of the next series examples. 
\begin{align*} 
\sum_{n \geq 0} \frac{(-1)^n}{(n^2+1)} & = 
     \frac{1}{2}\left(1 + \pi \csch(\pi)\right) \\ 
     & = 
     \sum_{j \geq 0} \frac{1}{2^{j+2}} \left(\binom{j+\imath}{\imath}^{-1} + 
     \binom{j-\imath}{-\imath}^{-1}\right) \\ \bigskip \\ 
\sum_{n \geq 0} \frac{(-1)^n}{(n^2+1)^2} & = 
     \frac{1}{4}\left(2 + \pi(1+\pi\coth(\pi)) \csch(\pi)\right) \\ 
     & = 
     \sum_{j \geq 0} \frac{1}{2^{j+3}}\left( 
     \binom{j+\imath}{\imath}^{-1} 
     \left(2 + \imath \HZNumber{j}{1}{1}{\imath}\right) - 
     \binom{j-\imath}{-\imath}^{-1} 
     \left(-2 + \imath \HZNumber{j}{1}{1}{-\imath}\right) 
     \right) \\ \bigskip \\ 
\sum_{n \geq 0} \frac{(-1)^n}{(n^2+1)^3} & = 
     \frac{1}{32}\left(16 + 6\pi (1+\pi\coth(\pi)) \csch(\pi) + 
     \pi^3 (3 + \cosh(2\pi)) \csch(\pi)^3\right) \\ 
     & = \phantom{+} 
     \sum_{j \geq 0} \frac{1}{2^{j+5}} \binom{j+\imath}{\imath}^{-1} \left( 
     8 + 5 \imath \HZNumber{j}{1}{1}{\imath} - \HZNumber{j}{1}{1}{\imath}^2 
     - \HZNumber{j}{2}{1}{\imath}\right) \\ 
     & \phantom{=\ } + 
     \sum_{j \geq 0} \frac{1}{2^{j+4}} \binom{j-\imath}{-\imath}^{-1} \left( 
     8 - 5 \imath \HZNumber{j}{1}{1}{-\imath} - \HZNumber{j}{1}{1}{-\imath}^2 
     - \HZNumber{j}{2}{1}{-\imath}\right) 
\end{align*} 
Another quadruple of trigonometric function series providing 
additional examples of 
expanding sums with quadratic denominators by partial fractions is 
given as follows \citep[\S 1.2]{SPECIAL-FUNCTIONS}: 
\begin{align*} 
\frac{\pi}{\sin(\pi x)} & = \frac{1}{x} + 
     \sum_{n \geq 0} \frac{2 (-1)^{n+1}}{x^2-(n+1)^2} \\ 
     & = \frac{1}{x} + \sum_{j \geq 0} \left[ 
     \binom{j+1-x}{1-x}^{-1} \frac{1}{1-x} - 
     \binom{j+1+x}{1+x}^{-1} \frac{1}{1+x}\right] \frac{1}{2^{j+1}} \\ 
     & = 
     \frac{1}{x} + \frac{\pi^2}{6} x + \frac{7 \pi^4}{360} x^3 + 
     \frac{31 \pi^6}{15120} x^5 + \frac{127 \pi^8}{604800} x^7 + O(x^9) \\ 
     \bigskip \\ 
\frac{1+x\csc(x)}{x} & = \sum_{n \geq 0} \frac{2x (-1)^n}{x^2-\pi^2 n^2} \\ 
     & = \sum_{j \geq 0} \left[\binom{j+\frac{x}{\pi}}{\frac{x}{\pi}}^{-1} 
     \frac{1}{x} + \binom{j-\frac{x}{\pi}}{-\frac{x}{\pi}}^{-1} \frac{1}{x} 
     \right] \frac{1}{2^{j+1}} \\ 
     & = 
     \frac{2}{x} + \frac{x}{6} + \frac{7}{360} x^3 + \frac{31}{15120} x^5 + 
     \frac{127}{604800} x^7 + O(x^9) \\ 
     \bigskip \\ 
\pi \sec(\pi x) & = \sum_{n=-\infty}^{\infty} \frac{(-1)^n}{n+x+\frac{1}{2}} \\ 
     & = 
     \frac{1}{\left(x+\frac{1}{2}\right)} - \sum_{j \geq 0} 
     \binom{j+x+\frac{3}{2}}{x + \frac{3}{2}}^{-1} 
     \frac{1}{2^{j+1} \left(x+\frac{3}{2}\right)} \\ 
     & \phantom{=\frac{1}{\left(1+\frac{1}{2}\right)}\ } + \sum_{j \geq 0} 
     \binom{j+\frac{1}{2}-x}{\frac{1}{2}-x}^{-1} 
     \frac{1}{2^{j+1} \left(\frac{1}{2}-x\right)} \\ 
     & = 
     \pi + \frac{\pi^3}{2} x^2 + \frac{5 \pi^5}{24} x^4 + 
     \frac{61 \pi^7}{720} x^6 + \frac{277 \pi^9}{8064} x^8 + O(x^{10}) \\ 
     \bigskip \\ 
\frac{2+x (1+x \cot(x)) \csc(x)}{4 x^4} & = 
     \sum_{n \geq 0} \frac{(-1)^n}{(x^2-\pi^2 n^2)^2} \\ 
     & = 
     \sum_{b = \pm 1} \sum_{n \geq 0} \frac{(-1)^n}{4}\left( 
     \frac{b}{x^3 (\pi n+bx)} + \frac{1}{x^2 (\pi n+bx)^2} 
     \right) \\ 
     & = 
     \sum_{b = \pm 1} \sum_{j \geq 0} \binom{j+\frac{bx}{\pi}}{j}^{-1} \left( 
     \frac{2}{x^4} + \frac{b}{x^3} \HZNumber{j}{1}{\pi}{bx} 
     \right) \frac{1}{2^{j+3}} \\ 
     & = 
     \frac{1}{x^4}-\frac{7}{720}-\frac{31 x^2}{15120} - 
     \frac{127 x^4}{403200}-\frac{73 x^6}{1710720}+O(x^7). 
\end{align*}

\begin{remark}[Expansions of General Zeta Series with Quadratic Denominators] 
For general quadratic zeta series of the form 
\begin{align} 
\label{eqn_GenQuadraticZetaSeries_def}
\sum_{n \geq 0} \frac{(-1)^n z^n}{(an^2+bn+c)^{s}}, 
\end{align} 
for integers $s \geq 1$ and constants $a, b, c \in \mathbb{R}$, we can 
apply the same procedure of factoring the denominators into linear factors of 
$z$ and taking partial fractions to write the first sums as a finite sum over 
functions of the form $\Phi(z, s, \alpha, \beta)$ for variable parameters, 
$\alpha, \beta \in \mathbb{C}$. 
We also note that a less standard definition of the 
Lerch transcendent function, $\Phi(z, s, a)$, is given by 
\begin{align*} 
\Phi^{\ast}(z, s, a) & = \sum_{n \geq 0} \frac{z^n}{\left[(z+a)^2\right]^{s/2}}, 
\end{align*} 
which also suggests an approach to a reduction to non--integer order 
exponents $s$ from the general quadratic zeta series in the forms of 
\eqref{eqn_GenQuadraticZetaSeries_def}. 
\end{remark} 

\subsection{Special Series Identities for the Riemann Zeta Function} 
\label{subSection_Examples_SeriesIdents_RZeta} 

\subsubsection{Series Generating the Riemann Zeta Function at the Even Integers} 

We first consider the following series identity for the 
zeta function constant, $\zeta(3)$, given by \citep[\S 7.10.2]{INCGAMMAFNS-APPS} 
\begin{align*}
\zeta(3) & = \frac{2 \pi^2}{9}\left(\Log(2) + 2 \times 
     \sum_{k \geq 0} \frac{\zeta(2k)}{2^{2k} (2k+3)} 
     \right) 
\end{align*} 
Since we know the next ordinary generating function as 
\citep[\S 25.8]{NISTHB} 
\begin{align*} 
\sum_{k \geq 0} \zeta(2k) z^{k} & = -\frac{\pi \sqrt{z}}{2} \cot(\pi \sqrt{z}), 
\end{align*} 
where $D_z^{(n)}[f(z) \cdot g(z)] = \sum_k \binom{n}{k} f^{(k)}(z) g^{(n-k)}(z)$, 
$D_z\left[\cot(z)\right] = -1 - \cot^2(z)$, and we can expand the $n^{th}$ 
derivatives of the cotangent function according to the known formula on the 
\emph{Wolfram Functions} website as 
\begin{align*} 
D_z^{(n)}[\cot(z)] & = \cot(z) \cdot \Iverson{n = 0} - 
     \csc^2(z) \cdot \Iverson{n = 1} \\ 
     & \phantom{=} - n \times 
     \sum_{0 \leq j < k < n} \binom{2k}{j} \binom{n-1}{k} 
     \frac{(-1)^k 2^{n-2k} (k-j)^{n-1}}{(k+1) \cdot \sin^{2k+2}(z)} 
     \sin\left(\frac{n\pi}{2} + 2(k-j) z\right), 
\end{align*} 
we may expand 
variants of the first sum for $\zeta(3)$ related to the zeta function 
constants, $\zeta(2k+1)$, for integers $k \geq 1$. 
In particular, we expand the first sum as follows: 
\begin{align*} 
\zeta(3) & = \frac{2 \pi^2}{9}\left(\Log(2) + 2 \times 
     \sum_{j \geq 0} \binom{j+\frac{3}{2}}{\frac{3}{2}}^{-1} 
     \frac{\left(-\frac{1}{4}\right)^j}{3 j!} 
     D_z^{(j)}\left[ 
     -\frac{\pi \sqrt{z}}{2} \cot\left(\pi \sqrt{z}\right) 
     \right] \Biggr\rvert_{z = \frac{1}{4}}
     \right). 
\end{align*} 
Related expansions of other series for the zeta function 
constants over the odd positive integers $s \geq 3$ include the 
next identity for integers $n \geq 1$ \citep[\S 7.10.2]{INCGAMMAFNS-APPS}. 
\begin{align*} 
\zeta(2n+1) & = \frac{2 (-1)^n (2\pi)^{2n}}{(2n-1) 2^{2n} + 1} \left( 
     \sum_{1 \leq k < n} \frac{(-1)^{k-1} k}{(2n-2k+1)!} 
     \frac{\zeta(2k+1)}{\pi^{2k}} + 
     \sum_{k \geq 0} \frac{(2k)!}{(2n+2k+1)!} \frac{\zeta(2k)}{2^{2k}} 
     \right) 
\end{align*} 
We then arrive at a BBP--type series for the constant, $\zeta(5)$, in the 
following form when 
$(a_1, a_2, a_3, a_4, a_5) = \left(\frac{1}{24}, -\frac{1}{6}, \frac{1}{4}, 
                                   -\frac{1}{6}, \frac{1}{24}\right)$: 
\begin{align*}
\zeta(5) & = \frac{16 \pi^2}{147} \zeta(3) + \frac{32 \pi^4}{49} \times 
     \sum_{j \geq 0} \sum_{1 \leq i \leq 5} 
     \binom{j+\frac{i}{2}}{\frac{i}{2}}^{-1} \frac{a_i}{i} 
     \frac{\left(-\frac{1}{4}\right)^j}{j!} D_z^{(j)}\left[ 
     -\frac{\pi \sqrt{z}}{2} \cot\left(\pi \sqrt{z}\right) 
     \right] \Biggr\rvert_{z = \frac{1}{4}}
\end{align*} 
Suppose that the ordinary generating function for the sequence, 
$\langle c_k^{\ast} \rangle$, is denoted by $C(z)$. Then provided that the 
function $C(z)$ has $j^{th}$ derivatives with respect to $z$ for 
$1 \leq j \leq 5$, we can transform this generating function into a 
generating function enumerating the next sequence enumerated in the form of 
\begin{align*}
 & \sum_{k \geq 0} \left(32 k^5 + 240 k^4 + 680 k^3 + 900 k^2 + 548 k + 120 
     \right) c_k^{\ast} z^k = \\ 
     & \phantom{=} 
     120 C(z) + 2400 z C^{\prime}(z) + 5100 z^2 C^{\prime\prime}(z) + 
     2920 z^3 C^{(3)}(z) + 560 z^4 C^{(4)}(z) + 32 z^5 C^{(5)}(z), 
\end{align*} 
which implies that 
\begin{align*} 
\widetilde{C}_{2k}(z) & := 
     \sum_{k \geq 0} (2k+1)(2k+2)(2k+3)(2k+4)(2k+5) \zeta(2k) z^k \\ 
     & = 
     \scriptstyle
     8 \pi ^6 z^3 \csc ^6\left(\pi  \sqrt{z}\right)+4 \pi ^4 z^2 \left(11 \pi ^2
     z \cot ^2\left(\pi  \sqrt{z}\right)-60 \pi  \sqrt{z} \cot \left(\pi 
     \sqrt{z}\right)+75\right) \csc ^4\left(\pi  \sqrt{z}\right) \\ 
     & \phantom{=\ } \scriptstyle + 
     4 \pi ^2 z
     \left(-30 \pi ^3 z^{3/2} \cot ^3\left(\pi  \sqrt{z}\right)+2 \pi ^4 z^2
     \cot ^4\left(\pi  \sqrt{z}\right)+150 \pi ^2 z \cot ^2\left(\pi 
     \sqrt{z}\right)-300 \pi  \sqrt{z} \cot \left(\pi 
     \sqrt{z}\right)+225\right) \csc ^2\left(\pi  \sqrt{z}\right) \\ 
     & \phantom{=\ } \scriptstyle - 
     360 \pi \sqrt{z} \cot \left(\pi  \sqrt{z}\right). 
\end{align*} 
From the generating function expansion in the previous equation, we have 
another ``\emph{coerced}'' variant of a degree--$2$ BBP--type formula 
of the following form when the corresponding coefficient sets are defined to be 
$(b_1, b_2, b_3, b_4, b_5) = \frac{1}{3456} \left(-25, -160, 0, 160, 25\right)$ and 
$(c_1, c_2, c_3, c_4, c_5) = \frac{1}{3456} \left(6, 96, 216, 96, 6\right)$: 
\begin{align*}
\zeta(5) & = \frac{16 \pi^2}{147} \zeta(3) \\ 
     & \phantom{=\ } + 
     \frac{32 \pi^4}{49} \times 
     \sum_{j \geq 0} \sum_{1 \leq i \leq 5} \left[
     \binom{j+\frac{i}{2}}{\frac{i}{2}}^{-1} \left( 
     \frac{b_i}{i} + \frac{c_i}{i^2} + \frac{c_i}{i} \HZNumber{j}{1}{2}{i} 
     \right) \right] 
     \frac{\left(-\frac{1}{4}\right)^j}{j!} 
     D_z^{(j)}\left[\widetilde{C}_{2k}(z)\right] \Biggr\rvert_{z = \frac{1}{4}}. 
\end{align*} 

\subsubsection{Zeta Function Constants Defined by the 
               Alternating Hurwitz Zeta Function} 

We conclude the applications in this section with several series for the 
\emph{alternating Hurwitz zeta function}, 
$\zeta^{\ast}(s, \alpha, \beta) + \beta^{-s} = \Phi(-1, s, \alpha, \beta)$. 
In particular, we see that \citep[\S 25.11(x)]{NISTHB} 
\begin{align*} 
\sum_{n \geq 0} (-1)^n \left[\frac{1}{(3n+1)^s} - \frac{1}{(3n+2)^s}\right] & = 
     6^{-s} \left[\zeta\left(s, \frac{1}{6}\right) - 
     \zeta\left(s, \frac{2}{3}\right) + 
     \zeta\left(s, \frac{5}{6}\right) - 
     \zeta\left(s, \frac{1}{3}\right) 
     \right] \\ 
     & = 
     6^{-s}\left(2^s-2\right)\left(3^s-1\right) \cdot \zeta(s). 
\end{align*} 
Examples of the last series identity for the Riemann zeta function, $\zeta(s)$, 
when $s = 2, 3, 4$ are given in the following equations in terms of the 
special cases of the transformation coefficients expanded in 
Table \ref{table_TableOfGenCoeffs_ab31} and 
Table \ref{table_TableOfGenCoeffs_ab32}: 
\begin{align*} 
\frac{2 \pi^2}{27} & = \sum_{i = 1, 2} \sum_{j \geq 0} 
     \binom{j+\frac{i}{3}}{\frac{i}{3}}^{-1} \left[ 
     \frac{1}{i^2} + \frac{1}{i} \HZNumber{j}{1}{3}{i}\right] 
     \frac{(-1)^{i+1}}{2^{j+1}} \\ 
\frac{13}{18} \zeta(3) & = \sum_{i = 1, 2} \sum_{j \geq 0} 
     \binom{j+\frac{i}{3}}{\frac{i}{3}}^{-1} \left[ 
     \frac{1}{i^3} + \frac{1}{i^2} \HZNumber{j}{1}{3}{i} + 
     \frac{1}{2 i} \left(\HZNumber{j}{1}{3}{i}^2 + \HZNumber{j}{2}{3}{i}\right) 
     \right] \frac{(-1)^{i+1}}{2^{j+1}} \\ 
\frac{7 \pi^4}{729} & = \sum_{i = 1, 2} \sum_{j \geq 0} 
     \binom{j+\frac{i}{3}}{\frac{i}{3}}^{-1} \Biggl[ 
     \frac{1}{i^4} + \frac{1}{i^3} \HZNumber{j}{1}{3}{i} + 
     \frac{1}{2 i^2} \left(\HZNumber{j}{1}{3}{i}^2 + 
     \HZNumber{j}{2}{3}{i}\right) \\ 
     & \phantom{= \sum_{i = 1, 2} \sum_{j \geq 0} \qquad\ } + 
     \frac{1}{6 i} \left(\HZNumber{j}{1}{3}{i}^3 + 2 \HZNumber{j}{1}{3}{i} 
     \HZNumber{j}{2}{3}{i} + 3 \HZNumber{j}{3}{3}{i}\right) 
     \Biggr] \frac{(-1)^{i+1}}{2^{j+1}}. 
\end{align*} 

\section{Conclusions and Final Examples} 
\label{Section_Concl} 

We have defined and proved special cases of a generalized generating function 
transform generating modified zeta functions and special zeta series. 
In Section \ref{subSection_GenHZStirlingNumbers} we 
connected the harmonic number expansions of generalized 
Stirling numbers of the first kind to partial sums of the modified 
Hurwitz zeta function defined by 
\eqref{eqn_GenModifiedHZetaFn_def} and \eqref{eqn_ModLTransFn_series_def} 
of Section \ref{Section_Intro}. 
The primary source of our new examples and applications 
is the generalization, or at least significant corollary, to the 
generating function transformations proved in \citep{GFTRANS2016}. 
Section \ref{subSection_Intro_Examples} of the introduction and 
Section \ref{Section_Examples} suggest many more 
important applications of these 
generalized forms of the generating function transformations explored in the 
first article. 

\subsection{Relations to Stirling Numbers of the Second Kind} 
\label{subSection_Concl_S2Analog}

The second geometric series transformation identity stated in 
\eqref{eqn_HZCfStar_sum_def_v1} effectively provides a combinatorial 
motivation for the known series for the Lerch transcendent function 
given by (See \citep{GFTRANS2016}) 
\begin{align*} 
\Phi(z, s, \alpha, \beta) & = \sum_{k \geq 0} \left(\frac{-z}{1-z}\right)^{k+1} 
     \sum_{0 \leq m \leq k} \binom{k}{m} 
     \frac{(-1)^{m+1}}{(\alpha m+\alpha+\beta)^{s}}. 
\end{align*} 
The even more general forms of the transformation coefficients defined by 
\eqref{eqn_HZCfStar_sum_def_v1} 
are considered to be Stirling numbers of the second kind 
``\emph{in reverse}'' in the sense that we have another related generalized 
form of the generating function transformation motivating the 
explorations in the first article defined by \citep[\cf \S 7.4]{GKP} 
\begin{align*}
\gkpSII{k}{j}_{\alpha,\beta} & := \sum_{0 \leq m \leq j} 
     \binom{j}{m} \frac{(-1)^{j-m} (\alpha m+\beta)^{k}}{j!} && \implies \\ 
\sum_{n \geq 0} (\alpha n+\beta)^{k} z^n & = 
     \sum_{0 \leq j \leq k} \gkpSII{k}{j}_{\alpha,\beta} 
     \frac{z^j \cdot j!}{(1-z)^{j+1}}. 
\end{align*} 
Moreover, we can extend this analog by observing that we also have the 
following negative-order identity involving the 
generalized Stirling numbers of the second kind defined by the last 
power series transformation identity \citep[\cf \S 26.8(v)]{NISTHB}: 
\begin{align*} 
\sum_{0 \leq j \leq n} (\alpha j+\beta)^{k} z^j & = 
     \sum_{0 \leq j \leq k} \gkpSII{k}{j}_{\alpha,\beta} z^j \times 
     D_z^{(j)}\left[\frac{1-z^{n+1}}{1-z}\right]. 
\end{align*} 

\subsection{Limitations and Versatility of the Transformations} 

\subsubsection{Some Limitations} 

Most of the generating functions we have employed in constructing the 
examples and applications within this article are based on variants of the 
geometric series, $G(z) = 1 / (1-cz)$, for some non-zero constant 
$c \in \mathbb{C}$ such that $|cz| < 1$ or when $cz \equiv -1$. 
One notable and obvious limitation of applying these geometric-series-based 
cases of our new transformations defined in 
Section \ref{subSection_Intro_GenGFSeriesEnums} is that 
we are not able to handle series of the form $\sum_{n} z^n / (\alpha n+\beta)^s$ 
when $|z| \equiv 1$, nor when $|\frac{1}{z}| < 1$. 
This restriction prevents us from constructing further examples of 
series for special zeta functions and constants such as 
\begin{align*} 
\frac{\pi^2}{8} & = \sum_{k \geq 0} \frac{1}{(2k+1)^2} = 
     1 + \frac{1}{3^2} + \frac{1}{5^2} + \frac{1}{7^2} + \cdots \\ 
\zeta(s) & = \frac{1}{1-2^{-s}} \times \sum_{n \geq 0} \frac{1}{(2n+1)^s}. 
\end{align*} 
On the other hand, expansions for multiples of $\pi^2$ and 
\emph{Catalan's constant}, $G$, for example, such as 
\begin{align*} 
\frac{\pi^2}{12} & = 
     \sum_{k \geq 0} \frac{(-1)^k}{(k+1)^2} \\ 
     & = 
     \sum_{k \geq 0} (-1)^k \left[ 
     \frac{13}{(3k+1)^2} - \frac{13}{(3k+2)^2} + \frac{4}{(3k+3)^2} 
     \right], 
\end{align*} 
are readily handled by our new transformations of ordinary power series 
generating functions. 

\subsubsection{Examples of Geometric-Series-Based 
               Generating Function Variants} 

Another example of a geometric-series-based zeta series variant is given by 
\begin{align*} 
\tan^{-1}(x) & = \sum_{n \geq 0} \frac{(-1)^n}{5^n} 
     \frac{F_{2n+1} t^{2n+1}}{(2n+1)} \\ 
     & = \frac{\sqrt{5}}{2 \imath} \times \sum_{b = \pm 1} \sum_{j \geq 0} 
     \binom{j+\frac{1}{2}}{j}^{-1} \frac{b}{\sqrt{5}} \left[ 
     \frac{\left(b\imath\varphi t / \sqrt{5}\right)^j}{ 
     \left(1-\frac{b\imath\varphi t}{\sqrt{5}}\right)^{j+1}} - 
     \frac{\left(b\imath\Phi t / \sqrt{5}\right)^j}{ 
     \left(1+\frac{b\imath\Phi t}{\sqrt{5}}\right)^{j+1}} 
     \right], 
\end{align*} 
for $t \equiv 2x / \left(1+\sqrt{1+\frac{4}{5} x^2}\right)$, where 
$F_{2n+1}$ denotes the \emph{$(2n+1)^{th}$ Fibonacci number} whose 
generating function is expanded in partial fractions as follows 
for $\varphi, \Phi := \frac{1}{2}\left(1 \pm \sqrt{5}\right)$ and the 
real-valued constants $c_1 := 1 / \sqrt{5}$ and $c_2 := -1 / \sqrt{5}$: 
\begin{align*} 
\sum_{n \geq 0} F_{2n+1} z^{2n+1} & = \frac{1}{2} \cdot \left( 
     \frac{c_1}{1-\varphi z} - \frac{c_1}{1+\varphi z} - 
     \frac{c_2}{1+\Phi z} + \frac{c_2}{1-\Phi z} 
     \right). 
\end{align*} 
We can also form yet other variants of the geometric-series-based 
transformations by considering Fourier series for special polynomials, 
such as the \emph{periodic Bernoulli polynomials}, 
$\widetilde{B}_n(x) \equiv B_n(x-\{x\})$, and the 
\emph{Euler polynomials}, $E_n(x) = n! \cdot [t^n] 2e^{tx} / (e^t+1)$, 
given in the forms of the following particular series expansions for 
$n \geq 1$ \citep[\S 24.8(i)]{NISTHB}: 
\begin{align*} 
\frac{E_{2n-1}(x)}{(2n-1)!} & = \frac{4 (-1)^n}{\pi^{2n}} \times 
     \sum_{k \geq 0} \frac{\cos\left((2k+1) \pi x\right)}{(2k+1)^{2n}} \\ 
\frac{E_{2n}(x)}{(2n)!} & = \frac{4 (-1)^n}{\pi^{2n+1}} \times 
     \sum_{k \geq 0} \frac{\sin\left((2k+1) \pi x\right)}{(2k+1)^{2n+1}}. 
\end{align*} 

\subsubsection{Series Involving Reciprocals of Binomial Coefficients} 

We can extend the method employed in constructing the results in 
Section \ref{subSection_Examples_SeriesIdents_RZeta} 
to further zeta series enumerating reciprocals of the 
\emph{central binomial coefficients}, $\binom{2n}{n}$, by first noticing that 
we have the following generating function: 
\begin{align*} 
\sum_{n \geq 0} \frac{z^{2n}}{\binom{2n}{n}} & = 
     4 \left(\frac{1}{4-z^2} + \frac{z}{(4-z^2)^{3/2}} 
     \sin^{-1}\left(\frac{z}{2}\right) 
     \right),\ |z| < 2. 
\end{align*} 
Though a formula for the $j^{th}$ derivatives of the right-hand-side function 
in the last equation is not clear, we may proceed to formulate the 
corresponding modified zeta series transformations involving the 
power series expansions of this ordinary generating function. 
Examples of applications of expanding series of this type 
through the new generating function transformations include the 
following three series \citep[\S 25.6(iii)]{NISTHB}: 
\begin{align*} 
\zeta(3) & = \frac{5}{2} \times \sum_{k \geq 1} 
     \frac{(-1)^{k-1}}{k^3 \binom{2k}{k}} \\ 
\zeta(2) - \csch^{-1}(2) \sinh^{-1}(2) & = 
     \sum_{n \geq 0} \frac{(-1)^n}{(2n+1)^2} \binom{2n}{n}^{-1} \\ 
\left(\sin^{-1}(x)\right)^2 & = \frac{1}{2} \times 
     \sum_{n \geq 0} \frac{(2x)^{2n}}{n^2 \binom{2n}{n}}. 
\end{align*} 

\subsubsection{Other Examples and Applications} 

One last class of applications of the modified zeta series transformations 
that is important to mention comprises geometric and exponential-series-based 
generating functions for which the parameters $(\alpha, \beta)$ in the 
coefficients from \eqref{eqn_HZCfStar_sum_def_v1} correspond to the 
expansion variable in a power series for a special function. 
We again demonstrate by example 
\citep[\cf \S 5.9(i)]{NISTHB}: 
\begin{align*} 
\Gamma(z) & = \sum_{n \geq 0} \frac{(-1)^n}{(n+z) \cdot n!} + 
     \int_1^{\infty} t^{z-1} e^{-t} dt \\ 
     & = 
     \sum_{j \geq 0} \frac{e^{-1}}{z} \binom{j+z}{z}^{-1} + 
     \int_1^{\infty} t^{z-1} e^{-t} dt \\ 
\sum_{k \geq 1} \frac{z}{(kz+1)^2} & = \sum_{k \geq 0} B_k z^k = 
     \int_0^{\infty} \frac{tz e^{-t}}{e^{tz}-1} dt \\ 
     & = 
     \sum_{k \geq 1} \frac{(-1)^k z}{(kz+1)^2} + 2 \times 
     \sum_{k \geq 0} \frac{z}{((2k+1) z+1)^2} \\ 
     & = 
     \sum_{k \geq 1} \frac{(-1)^k z}{(kz+1)^2} + 2 \times 
     \sum_{k \geq 0} \frac{(-1)^k z}{((2k+1) z+1)^2} + 4 \times 
     \sum_{k \geq 0} \frac{z}{((4k+3) z+1)^2} \\ 
     & = 
     \sum_{i \geq 0} \sum_{k \geq 0} \frac{2^i (-1)^k z}{ 
     \left((2^i k + 2^i-1) z + 1\right)^2}. 
\end{align*} 

\subsection{For Readers and Reviewers} 

A summary \emph{Mathematica} notebook providing numerical data and 
supporting computations in deriving key results and new applications to 
specific series is provided online at the following 
\emph{Google Drive} link: 
\url{https://drive.google.com/file/d/0B6na6iIT7ICZMjJnOFcySmlBMGs/view?usp=sharing}. 
The intention of this supplementary document included with the submission of 
this article is to help the reviewer process the article more quickly, and 
to assist the reader with verifying and modifying the examples 
presented as applications of the new results cited above. 


\renewcommand{\refname}{References}

\end{document}